\documentclass{svproc}

\usepackage{url}
\usepackage{xcolor} 
\usepackage{amsmath}
\usepackage{xfrac}
\usepackage{ulem} 

\usepackage{soul}

\begin{document}
\newtheorem{obs}[theorem]{Observation}
\newtheorem{conj}[theorem]{Conjecture}
\newtheorem{prop}[theorem]{Proposition}
\newtheorem{lem}[theorem]{Lemma}
\newtheorem{defn}[theorem]{Definition}
\newtheorem{cor}[theorem]{Corollary}
\newtheorem{ques}[theorem]{Question}
\newtheorem{prob}[theorem]{Problem}

\mainmatter              

\title{Constructions of betweenness-uniform graphs from trees}
\titlerunning{}  
%
\author{David Hartman\inst{1,2} \and Aneta Pokorná\inst{1}}
\authorrunning{Hartman, Pokorná} 
%
\tocauthor{David Hartman and Aneta Pokorná}
\institute{
Computer Science Institute of Charles University, Faculty of Mathematics and Physics, Charles University, 118 00, Prague, Czech Republic\\
\url{https://iuuk.mff.cuni.cz}
\and
The Institute of Computer Science of the Czech Academy of Sciences, Prague, Czech Republic
}

\maketitle              

\begin{abstract}
Betweenness centrality is a measure of the importance of a vertex $x$ inside a network based on the fraction of shortest paths passing through $x$. 
We study a blow-up construction that has been shown to produce graphs with uniform distribution of betweenness. 
We disprove the conjecture about this procedure’s universality by showing that trees with a diameter at least three cannot be transformed into betweenness-uniform by the blow-up construction.
It remains open to characterize graphs for which the blow-up construction can produce betweenness-uniform graphs.


\end{abstract}
In our vibrant society, everything is moving. Goods are being transported between factories, warehouses and stores, ideas are communicated among people, data are passing through the Internet. Such
transfers are often realized via shortest paths, and thus the structure of the underlying network plays an essential part in the workload of particular nodes. One way to measure such expected workload of a node in a network is betweenness centrality. This measure is based on the fraction of shortest path passing through a given vertex.
More precisely, for a graph $G = (V,E)$ with $x \in V(G)$ we define {\it betweenness centrality} as
$$B(x) = \sum_{\{u,v\} \in {V(G)\setminus \{x\} \choose 2}} \frac{\sigma_{u,v}(x)}{\sigma_{u,v}},$$
where $\sigma_{u,v}$ denotes the number of shortest paths between $u$ and $v$ and $\sigma_{u,v}(x)$ is the number of shortest paths between $u$ and $v$ passing through $x$.\cite{freeman}
Indeed, it seems to be a highly useful notion of measuring the importance
of nodes, as shown by its numerous applications in neuroscience\cite{brain}, chemistry\cite{chemistry}, sociology\cite{social}, and transportation\cite{temporal-bc}.
From a more theoretical perspective, betweenness centrality has also been shown to be a helpful measure in modelling random planar graphs\cite{bc-rand-planar}.

Once we measure the tendency to put uneven workload to different vertices, we might want to optimize the underlying network such that the communication is spread more evenly. Apart from preventing overload and potential collapse, it might also be important from the strategic perspective, as we do not want to have a single point of failure.
The extremal case is then the class of graphs with uniform distribution of betweenness.
A graph is called {\it betweenness-uniform}, if the value of betweenness is the same for all vertices.
It is easily observable that vertex-transitive graphs are betweenness-uniform.\cite{Pokorna2020}
However, for any fixed $n$, there are superpolynomially many betweenness-uniform graphs, which are not vertex-transitive\cite{GagoCoronicovaHurajovaMadaras2013}.
Apart from the fact that the class of distance-regular graphs is betweenness-uniform \cite{qgt,GagoCoronicovaHurajovaMadaras2013}, not much more is known about the characterization of betweenness-uniform graphs.
An important property of betweenness-uniform graphs is that they are always $2$-connected.\cite{GagoCoronicovaHurajovaMadaras2013} 
Actually, any betweenness-uniform graph is $3$-connected, unless it is isomorphic to a cycle.\cite{caldam2021}
There exist other studies concerned with both the values of betweenness in betweenness-uniform graphs\cite{bu-values} and edge betweenness-uniform graphs\cite{edge-bug}.

Throughout this text we use a standard notation.
A graph $G$ has a vertex set $V(G)$ and edge set $E(G)$.
The {\it degree} of a vertex $v \in V(G)$ is $\deg_G(v) = |\{u : uv \in E(G)\}| = |N_G(v)|$. Subscript is omitted whenever $G$ is clear from the context.
The {\it distance} $d(u,v)$ between $u,v \in V(G)$ is the length of the shortest path connecting these vertices.
We say that a graph $G$ has {\it diameter} $k$ if $k = \max_{u,v \in V(G)} \{d(u,v)\}$.
We denote by $P_n$ the path on $n$ vertices, by $K_n$ the complete graph on $n$ vertices and by $I_n$ the edge-less graph on $n$ vertices
and by $S_{n-1}$ the star with one central vertex and $n-1$ leafs.
We write $G \cong H$ whenever graph $G$ is isomorphic to graph $H$.
The set $\{1, \dots, n\}$ is denoted by $[n]$.

Let $G$ be a graph with vertex set $\{v_1, \dots, v_n\}$ and $H_1, \dots, H_n$ 
be a set of graphs. The graph $G[H_1, \dots, H_n]$ is 
defined on a vertex set
obtained from $V(G)$ via substituting each $v_i \in V(G)$ by
set of vertices $V(H_i)$ and has edge set $E(G[H_1, \dots, H_n]) = \bigcup_{i=1}^n E(H_i) \cup \{uv\,|\,u \in V(H_i), v \in V(H_j), v_iv_j \in E(G)\}$.
We will call $H_i$ as the {\it blow-up} of $v_i$.

Our primary motivation in this work is to study graphs potentially
generated via the following conjecture.
\begin{conj}[Coroničová Hurajová, Gago, Madaras, 2013 \cite{GagoCoronicovaHurajovaMadaras2013}]
For any graph $G$ on a vertex set $\{v_1, \dots , v_n\}$ there exist graphs
$H_1, \dots, H_n$ such that $G[H_1 , \dots , H_n]$ is betweenness-uniform.
\end{conj}
In the following, we deal with generating graphs from paths. 
We show the new blow-up constructions for $P_3$ and stars leading to betweenness-uniform graphs.
Contrary to that, we show that no such construction exists for $P_4$
and further generalize this claim to trees with diameter at least three,
disproving thus the conjecture.

By $v \in V(H_i)$ we mean that $v$ is also in a subgraph of $G[H_1, \dots, H_n]$
defined on vertices blown-up from $v_i$.
Whenever we consider a graph $G[H_1, \dots, H_n]$, we assume that the 
underlying graph $G$ has at least two vertices and is connected.
From these assumptions it is not hard to realize that the distance of any two vertices inside the same blow-up graph $H_i$ is at most two.

While considering betweenness of a vertex $v$ in $G[H_1, \dots, H_n]$, we will distinguish the paths contributing to $B(v)$ according to their endpoints
 $x, y$. Either $x,y \in V(H_i)$ for some $i \in [n]$, or $x \in V(H_i), y \in V(H_j)$ for $i,j \in [n]$ such that $i \neq j$.
We say that the first type of paths contributes to the {\it local} betweenness of $v$ for given $H_i$, $B^{H_i}(v)$, and that the second type of paths contributes to the ${\it global}$ betweenness of $v$, $B^G(v)$.
\begin{obs}\label{o:bc-expr}
For $v \in V(H_i)$, 
$B(v) = B^G(v) + B^{H_i}(v) + \sum_{j: v_j \in N_G(v_i)} B^{H_j}(v).$
\end{obs}

This observation enables us to ignore any paths longer than two between $x,y \in V(H_i)$ for any $i$, as such paths do not influence betweenness of any vertex. 
As a result, the only long paths influencing betweenness of any vertex are the paths between distinct $H_i$'s. 
Furthermore, the shortest paths between $H_i$ and $H_j$ for $j \neq i$ contain at most one vertex from each $H_k$, $k \in [n]$.

Our goal is to find a betweenness uniform graph using the blow-up construction for a path $P_k$. 
We iterate candidates via different values of $k$.
Starting from the smallest paths, $P_2 = K_2$ is betweenness-uniform and thus
$P_2[K_m, K_m]$ is betweenness-uniform for any integer $m$
as has been shown in the article \cite{GagoCoronicovaHurajovaMadaras2013}.

Even though the path $P_3$ and star $S_k$ are not betweenness-uniform, but we still can find a blow-up resulting in a betweenness-uniform graph.
\begin{obs}
For $P_3$, the graph $P_3[K_1, I_2, K_1]$ is isomorphic to $C_4$ and thus betweenness-uniform. 
Furthermore, $P_3[I_a, I_{a+b}, I_b]$ is betweenness-uniform for any $a,b$ positive integers.
\end{obs}

\begin{obs}
For any star $S_k$ with vertices $\{v_1, \dots, v_k\}$ of degree one and $\deg(v_{k+1}) = k$ there exists a blow-up construction 
$S_k[I_{s_1}, \dots, I_{s_k}, I_{\sum_{j=1}^k s_j}]$ giving a betweenness-uniform graph.
\end{obs}

Before moving to longer paths, we observe some properties of the blow-up graphs. 
From the fact that any betweenness-uniform graph is $2$-connected\cite{GagoCoronicovaHurajovaMadaras2013}, $|V(H_i)| > 1$ for any $v_i \in V(G)$ of degree at least two.

Let $u,v \in V(G[H_1, \dots, H_n])$ such that $u,v \in V(H_i)$ and $G' = G[H_1, \dots, H_n]$. We denote $\sigma_{u,v}^{H_i, G'}$ the number of $uv$-paths inside $H_i$, which have the same length as the shortest $uv$-path in $G'$. Specifically, $\sigma_{u,v}^{H_i, G'}$ is one if $uv \in E(H_i)$ and zero whenever $uv \notin E(H_i)$
and there is no $uv$-path of length two in $H_i$.

For simplicity, we denote $n_i = \sum_{j: v_j \in N_G(v_i)} |V(H_j)|$ the sum of sizes of the neighbouring blow-up graphs.

\begin{obs}\label{c:BH-size}
Let $G' = G[H_1, \dots, H_n]$, $x \in V(H_i)$ for some $i \in [n]$ and $v_iv_j \in E(G)$. Then 
$$B^{H_j}(x) = \sum_{\substack{u,v \in V(H_j)\\uv \notin E(H_j)}} \frac{1}{\sigma_{u,v}^{H_j, G'} + \sum_{\substack{k:\\ v_k \in N_G(v_j)}} |V(H_k)|}
= \sum_{\substack{u,v \in V(H_j)\\uv \notin E(H_j)}} \frac{1}{\sigma_{u,v}^{H_j, G'} + n_j}
$$
\end{obs}

Now we focus on $G = P_4$ and denote $G' = P_4[H_1, H_2, H_3, H_4]$.
We start by expressing betweenness centrality of a vertex $x \in V(H_1)$ and $y \in V(H_2)$ for $\deg_{G}(v_1) = 1$ and $v_1v_2 \in E(G)$.

Using Observation~\ref{o:bc-expr}, $B(x) = B^G(x) + B^{H_1}(x) + B^{H_2}(x)$.
From $v_1$ being an endpoint we have $B^G(v_1) = 0$.
Furthermore, for $N_1(x) = N(x) \cap V(H_1)$,
$$B^{H_1}(x) = \sum_{\substack{x_1,x_2 \in N_1(x)\\ x_1x_2 \notin E(G')}} \frac{1}{\sigma_{x_1x_2}^{H_1,G'} + |V(H_2)|}
\text{ and }
B^{H_2}(x) = \sum_{\substack{y_1y_2 \in V(H_2)\\ y_1y_2 \notin E(G')}} \frac{1}{\sigma_{y_1y_2}^{H_2, G'} + n_2}.$$
There are no other paths contributing to $B(x)$ for $x \in V(H_1)$.

We can use Observation~\ref{o:bc-expr} and Observation~\ref{c:BH-size} to express the betweenness of the vertex $y$,
$$B(y) = B^G(y) + B^{H_2}(y) + \sum_{j: v_j \in N_G(v_2)} \sum_{\substack{u,v \in V(H_j),\\ uv \notin E(H_j)}} 
\frac{1}{\sigma_{u,v}^{H_j, G'} + n_j}$$
The global betweenness of $y$ is closely related to the betweenness of $v_2$ in $G$.
All paths between $V(H_1)$ and $V(G') \setminus \{V(H_1) \cup V(H_2)\}$
contribute to each vertex of $V(H_2)$ by the same amount, so 
$B^G(y) = \frac{|V(H_1)| \cdot \big(|V(G')| - |V(H_1)| - |V(H_2)|\big)}{|V(H_2)|}.$
An interesting conclusion is that the global betweenness of $y$ grows heavily with the growing size of $|V(G')| - |V(H_2)|$ and gets small when $V(H_2)$ forms a large fraction of $V(G')$. 

The local contributions to $y$ behave similarly to the local contributions of $x$:
$$B^{H_2}(y) = \sum_{\substack{y_1y_2 \in N(y) \cap V(H_2)\\y_1y_2 \notin E(G')}} \frac{1}{\sigma_{y_1y_2}^{H_2, G'} + n_2}
\text { and }
B^{H_j}(y) = \sum_{\substack{y_1y_2 \in V(H_j)\\ y_1y_2 \notin E(G')}} \frac{1}{\sigma_{y_1y_2}^{H_j, G'} + n_j}$$
for each  $j$ such that $v_j \in N_G(v_2)$.

Note that the sizes of the local contributions shrink with the growing sizes of $H_i$'s in their neighbourhood, because the denominator grows.
\begin{obs}
For $deg_G(v_1) = 1$ and $v_1v_2 \in E(G)$, $B^{H_1}(x) \leq B^{H_1}(y)$.
\end{obs}

Using previous observation and the definition of $B^{H_i}$ we infer
\begin{equation}
B^{H_1}_{y-x} = B^{H_1}(y) - B^{H_1}(x) = \sum_{\substack{y_1 \in V(H_1), y_2 \in V(H_1)\setminus N(x)\\y_1y_2 \notin E(G')}} \frac{1}{\sigma_{y_1y_2}^{H_1, G'} + |V(H_2)|}.\label{e:BH1yx}
\end{equation}
Similarly, we obtain $B^{H_2}(x) \geq B^{H_2}(y)$ and thus 
\begin{equation}
B^{H_2}_{x-y} = B^{H_2}(x) - B^{H_2}(y) = \sum_{\substack{y_1 \in V(H_2), y_2 \in V(H_2)\setminus N(y)\\y_1y_2 \notin E(G')}} \frac{1}{\sigma_{y_1y_2}^{H_2, G'} + n_2}.\label{e:BH2xy}
\end{equation}

When considering $B(x) = B(y)$ and using equations (\ref{e:BH1yx}) and (\ref{e:BH2xy}) we obtain
\begin{align}
    B^{H_2}_{x-y} =&\ B^G(y) + B^{H_1}_{y-x} + \sum_{j: v_j \in 
    N_G(v_2), j\neq 1} B^{H_j}(y) \label{e:BxBy}
\end{align}
We will show that even when trying to maximize $B(x)$, it will always be smaller than $B(y)$ for sufficiently large $G'$.

We denote $$\Delta(x,y) = \frac{B^{H_2}_{x-y}}{B^G(y) + B^{H_1}_{y-x} + \sum_{\substack{j: v_j \in N_G(v_2)\\j\neq 1}} B^{H_j}(y)}.$$
If $B(x) = B(y)$, $\Delta(x,y) = 1$. For $B(x) < B(y)$ we obtain $\Delta(x,y) < 1$ and vice versa.

\begin{lem}\label{o:leaf-neighbour-independent}
Let $m = |H_2|$ be fixed. Then $\Delta(x,y)$ is maximized for $H_2 \cong I_m$.
\end{lem}

\begin{lem}\label{o:leaf-clique}
Let $m = |H_1|$ be fixed. Then $\Delta(x,y)$ is maximized for $H_1 \cong K_m$.
\end{lem}

In order to produce betweenness uniform graph using the above-described blow-up construction we need to maximize $\Delta(x,y)$, i.e. maximize corresponding $B(x)$ and minimize $B(y)$.
From the lemmas above we obtain
$H_1 \cong K_a, H_2 \cong I_b$, $H_3 \cong I_c, H_4 \cong K_d$ for some 
positive integers $a,b,c$ and $d$.
By this assumption, equation~(\ref{e:BxBy}) simplifies. 
We transform it into an inequality, as by maximizing $\Delta(x,y)$ we must allow the case when $B(x) > B(y)$. 
Note that by adding edges into $H_3$ betweenness of vertices in both $H_2$ and $H_4$ decreases.
By using the inequality for both $x \in V(H_1), y \in V(H_2)$ and $x' \in V(H_4), y' \in V(H_3)$ we obtain
$$ \frac{{b \choose 2}}{(a+c)} \geq \frac{a(c+d)}{b} + 0 + \frac{{c \choose 2}}{(b+d)}
\text{ and }
\frac{{c \choose 2}}{(b+d)} \geq \frac{d(a+b)}{c} + 0 + \frac{{b \choose 2}}{(a+c)}$$
By substituting the second inequality to the first one we get
$$ 0 \geq \frac{a(c+d)}{b} + \frac{d(a+b)}{c} 
\text{ implying }
0 \geq ac(c+d) + bd(a+b)$$
which cannot be fulfilled by $a,b,c,d$ positive integers and thus 
$B(x) < B(y)$ and $G'$ is not betweenness-uniform.

By realizing that the global betweenness of $y$ grows with growing number of vertices in the path while the betweenness of the endpoints of the path stays the same, we obtain a result summarizing betweenness-uniform blow-up constructions of paths.
\begin{prop}\label{c:no-blowup-Pk}
For $P_k$ with $k \geq 4$ there is no blow-up construction resulting in a betweenness-uniform graph.
\end{prop}

\begin{theorem}
Let $G$ be a tree with diameter $d$ at least three and $|V(G)| = n$.
Then there are no graphs $H_1, \dots, H_n$ such that 
$G[H_1, \dots, H_n]$ is a betweenness-uniform graph.
\end{theorem}
The idea of the proof is that $P_d$ is contained is $G$ as a subgraph such that its endpoints have degree one. It is not hard to see that 
by adding subtrees adjacent to the internal vertices of the path, the betweenness of end-vertices is not increasing and betweenness of the internal vertices is not decreasing.

Our results show that there are some non-betweenness-uniform graphs, such as $P_3$ and $S_k$, which can be transformed to a betweenness-uniform graph by the blow-up construction, but there are other non-betweenness-uniform graphs, which cannot be transformed into a betweenness-uniform graph by the blow-up construction for any $H_1, \dots, H_n$.
\begin{prob}
Determine a class $\mathcal{B}$ of graphs that can be transformed into betweenness uniform via blow-up construction with a suitable choice of $H_1, \ldots, H_n$.
\end{prob}

\begin{conj}
Let $G$ be a graph of order $n$ with diameter at least three having a vertex cut of size one. Then there are no graphs $H_1, \dots, H_n$
such that $G[H_1, \dots, H_n]$ is betweenness-uniform.
\end{conj}

\subsubsection{Acknowledgments.} David Hartman and Aneta Pokorná were partially supported by ERC Synergy grant DYNASNET grant agreement no. 810115.


\bibliographystyle{splncs03}
\typeout{}
\bibliography{betweenness-blow-up}

\begin{thebibliography}{10}
\providecommand{\url}[1]{\texttt{#1}}
\providecommand{\urlprefix}{URL }

\bibitem{freeman}
Freeman, L.C.: A set of measures of centrality based on betweenness. Sociometry
   40,  35--41 (1977)

\bibitem{qgt}
Gago, S., Coroničová~Hurajová, J., Madaras, T.: Betweenness centrality in
  graphs. In: Dehmer, M., Emmert-Streib, F. (eds.) Quantitative Graph Theory:
  Mathematical Foundations and Applications. Discrete Mathematics and Its
  Applications, Taylor \& Francis (2014)

\bibitem{GagoCoronicovaHurajovaMadaras2013}
Gago, S., Hurajová-Coroničová, J., Madaras, T.: {On betweenness-uniform
  graphs}. {Czechoslovak Mathematical Journal}  {63}({3}),  {629--642} ({2013})

\bibitem{social}
Girvan, M., Newman, M.E.J.: Community structure in social and biological
  networks. Proceedings of the National Academy of Sciences  99(12),
  7821--7826 (2002)

\bibitem{caldam2021}
Hartman, D., Pokorn{\'a}, A., Valtr, P.: On the connectivity and the diameter
  of betweenness-uniform graphs. In: Mudgal, A., Subramanian, C.R. (eds.)
  Algorithms and Discrete Applied Mathematics. pp. 317--330. Springer, Cham
  (2021)

\bibitem{bu-values}
Joseph, S., Ajitha, V., Jose, B.K.: Construction of some betweenness uniform
  graphs. AIP Conference Proceedings  2261(1),  030015 (2020)

\bibitem{bc-rand-planar}
Kirkley, A., Barbosa, H., Barthelemy, M., Ghoshal, G.: From the betweenness
  centrality in street networks to structural invariants in random planar
  graphs. Nature Communications  9,  2041--1723 (2018)

\bibitem{edge-bug}
Newman, H., Miranda, H., Florez, R., Narayan, D.A.: Uniform edge betweenness
  centrality. EJGTA  8(2),  265–300 (2020)

\bibitem{brain}
Pacheco, L.: An fmri study on betweenness centrality in brain networks of
  teenagers with inhaled substance abuse disorder. Am J Biomed Sci \& Res
  4(3),  136--139 (2019)

\bibitem{Pokorna2020}
Pokorná, A.: Characteristics of network centralities. Master's thesis, Charles
  University, Prague, Czech Republic (2020)

\bibitem{temporal-bc}
Zaoli, S., Mazzarisi, P., Lillo, F.: Betweenness centrality for temporal
  multiplexes. Scientific Reports  11,  2045--2322 (2021)

\bibitem{chemistry}
Zhao, P., Nackman, S.M., Law, C.K.: On the application of betweenness
  centrality in chemical network analysis: Computational diagnostics and model
  reduction. Combustion and Flame  162(8),  2991--2998 (2015)

\end{thebibliography}

\end{document}